\documentclass{amsart}
\usepackage{geometry} 
\usepackage{url}
\usepackage{amsthm}
\usepackage{amssymb}
\usepackage{latexsym}
\usepackage{amsfonts}
\usepackage{amsmath}
\usepackage{amstext}
\usepackage{accents}
\usepackage{cancel}
\usepackage{eucal}
\usepackage{amscd}
\usepackage{hyperref}
\usepackage[shortalphabetic]{amsrefs}
\usepackage{tikz}
\usepackage[pagewise]{lineno}

\title[Curve diffusion of planar curves inside cones]{Curvature diffusion of planar curves with generalised Neumann boundary conditions inside cones}

\author{Mashniah A. Gazwani and James McCoy}

\date{27 December 2023} 
\begin{document}

\begin{abstract}
We study families of smooth immersed  regular planar curves    $ \alpha : \left [-1,1  \right ]\times \left [0,T  \right )\to \mathbb{R}^{2}$ satisfying the fourth order nonlinear curve diffusion flow with generalised Neumann  boundary conditions inside cones.  We show that if the initial curve has sufficiently small oscillation of curvature then this remains so under the flow.  Such families of evolving curves either exist for a finite time, when an end of the curve has reached the tip of the cone or the curvature has become unbounded in $L^2$, or they exist for all time and converge exponentially in the $C^{\infty}$- topology to a circular arc that, together with the cone boundary encloses the same area as that of the initial curve and cone boundary.  The same kind of result is possible for the higher order polyharmonic curve diffusion flows with appropriate boundary conditions; in particular in the sixth order case the smallness condition on the oscillation of curvature is exactly the same as for curve diffusion.\\  \\
 \textbf{keywords}: Curve diffusion flow,  fourth order geometric evolution equation, Neumann  boundary condition.
 \small{The research of the second author was supported by Discovery Project DP180100431 of the Australian Research Council.  The research of the first author was supported by a postgraduate scholarship from the Department of Mathematics, College of Science, Imam Abdulrahman Bin Faisal University, P. O. Box 1982, Dammam, Saudi Arabia.  The authors are grateful for this support.  The authors are also grateful to G. Wheeler and V. Wheeler for useful discussions and to B. Andrews for his interest in this work.  The authors express their gratitude to the anonymous referee for his or her comments that have led to improvements to the article.

\url{Mashniah.Gazwani@uon.edu.au, , magazwani@iau.edu.sa, James.McCoy@newcastle.edu.au}
 }
\end{abstract}
\maketitle
\section{Introduction} \label{S:intro}

Study of the curve diffusion flow of planar curves goes back at least to Mullins \cite{M57}, where it was proposed with the practical application of thermal grooving of a heated polycrystal.  Since then other applications have been proposed and the flow has been studied in various mathematical settings.  In particular, we refer to the analyses of closed planar curves evolving by curve diffusion flow by G. Wheeler \cite{W13} and curves with generalised Neumann boundary conditions between parallel lines by G. Wheeler and V. Wheeler \cite{WW} and the references contained therein.

It is natural to question in which other settings with boundary an analysis similar to that in \cite{WW} may yield a classification result for long-time behaviour of solutions to the curve diffusion flow.  In this article we consider a similar setup, however instead of parallel supporting lines we consider an initial curve between two rays that form a cone in the plane.  One reason for considering this setting is that, just as straight lines are stationary solutions for the curve diffusion flow between parallel lines with generalised Neumann boundary conditions, in cones it is clear that circular arcs are stationary solutions to the curve diffusion flow again with generalised Neumann boundary conditions.  In each setting, as with geometric flows more generally, one hopes solutions might converge in some sense to stationary solutions, at least if they are initially `close enough'.  The second reason is that, while we no longer have the obvious lower bound on the length of the curve as the distance between the two parallel supporting lines, we do know instead that the area bounded by the cone and the evolving curve remains constant under the flow (Corollary \ref{T:evlneqns2} (ii)) and this facilitates a crucial lower length bound under the flow.  On the other hand, our setting introduces additional complexity: we must rule out flow to the tip of the cone (at which time the problem ceases to be well-posed).  Moreover, since the average curvature of the evolving curve is no longer identically equal to zero, additional terms arise in the relevant evolution equations and these must be appropriately handled in the analysis.  Nevertheless despite these issues we are able to again prove a satisfactory classification in this setting, Theorems \ref{T:main1} and \ref{T:main2} below.  In particular, given an initial curve of sufficiently small oscillation of curvature sufficiently far from the tip of the cone, a unique solution to the curve diffusion flow with generalised Neumann boundary conditions exists for all time and converges exponentially in the $C^\infty$ topology to a unique circular arc.

\begin{figure}
    \centering
\tikzset{every picture/.style={line width=0.75pt}}
\begin{tikzpicture}[x=0.75pt,y=0.75pt,yscale=-1,xscale=1]
\draw    (98,84) -- (578.33,84) ;
\draw [color={rgb, 255:red, 80; green, 227; blue, 194 }  ,draw opacity=1 ]   (108.33,84) .. controls (96.33,55) and (226.33,55) .. (330.63,57.34) .. controls (434.92,59.69) and (584.33,50) .. (569.33,85) ;
\draw  [color={rgb, 255:red, 80; green, 227; blue, 194 }  ,draw opacity=1 ] (108.11,77.72) -- (115.89,77.45) -- (116.11,83.72) ;
\draw  [color={rgb, 255:red, 80; green, 227; blue, 194 }  ,draw opacity=1 ] (561.37,85.05) -- (561.33,78.33) -- (569.29,78.28) ;
\draw [color={rgb, 255:red, 208; green, 2; blue, 27 }  ,draw opacity=1 ]   (103,89) -- (129.33,89.93) ;
\draw [shift={(131.33,90)}, rotate = 182.02] [color={rgb, 255:red, 208; green, 2; blue, 27 }  ,draw opacity=1 ][line width=0.75]    (10.93,-3.29) .. controls (6.95,-1.4) and (3.31,-0.3) .. (0,0) .. controls (3.31,0.3) and (6.95,1.4) .. (10.93,3.29)   ;
\draw [color={rgb, 255:red, 208; green, 2; blue, 27 }  ,draw opacity=1 ]   (568.33,94.33) -- (536.33,94.02) ;
\draw [shift={(534.33,94)}, rotate = 0.56] [color={rgb, 255:red, 208; green, 2; blue, 27 }  ,draw opacity=1 ][line width=0.75]    (10.93,-3.29) .. controls (6.95,-1.4) and (3.31,-0.3) .. (0,0) .. controls (3.31,0.3) and (6.95,1.4) .. (10.93,3.29)   ;
\draw    (520.33,122) -- (520.33,99) ;
\draw [shift={(520.33,97)}, rotate = 90] [color={rgb, 255:red, 0; green, 0; blue, 0 }  ][line width=0.75]    (10.93,-3.29) .. controls (6.95,-1.4) and (3.31,-0.3) .. (0,0) .. controls (3.31,0.3) and (6.95,1.4) .. (10.93,3.29)   ;
\draw    (72.33,262) -- (540.33,259) ;
\draw  [color={rgb, 255:red, 80; green, 227; blue, 194 }  ,draw opacity=1 ] (101.54,252.74) -- (110.54,252.74) -- (110.54,260.07) ;
\draw [color={rgb, 255:red, 208; green, 2; blue, 27 }  ,draw opacity=1 ]   (79,270.33) -- (106.33,270.33) ;
\draw [shift={(108.33,270.33)}, rotate = 180] [color={rgb, 255:red, 208; green, 2; blue, 27 }  ,draw opacity=1 ][line width=0.75]    (10.93,-3.29) .. controls (6.95,-1.4) and (3.31,-0.3) .. (0,0) .. controls (3.31,0.3) and (6.95,1.4) .. (10.93,3.29)   ;
\draw    (541,299.33) -- (541.31,273.33) ;
\draw [shift={(541.33,271.33)}, rotate = 90.68] [color={rgb, 255:red, 0; green, 0; blue, 0 }  ][line width=0.75]    (10.93,-3.29) .. controls (6.95,-1.4) and (3.31,-0.3) .. (0,0) .. controls (3.31,0.3) and (6.95,1.4) .. (10.93,3.29)   ;
\draw [color={rgb, 255:red, 208; green, 2; blue, 27 }  ,draw opacity=1 ]   (577.33,269.33) -- (553.33,269.33) ;
\draw [shift={(551.33,269.33)}, rotate = 360] [color={rgb, 255:red, 208; green, 2; blue, 27 }  ,draw opacity=1 ][line width=0.75]    (10.93,-3.29) .. controls (6.95,-1.4) and (3.31,-0.3) .. (0,0) .. controls (3.31,0.3) and (6.95,1.4) .. (10.93,3.29)   ;
\draw    (560.33,262.33) ;
\draw  [line width=3] [line join = round][line cap = round] (520.33,83) .. controls (520.61,83.56) and (521.04,84.33) .. (521.67,84.33) ;
\draw  [line width=3] [line join = round][line cap = round] (521,85) .. controls (521.94,85) and (519.91,83.23) .. (519,83) .. controls (517.8,82.7) and (520.33,86.41) .. (520.33,84.33) ;
\draw  [color={rgb, 255:red, 80; green, 227; blue, 194 }  ,draw opacity=1 ] (559.41,254.5) -- (556.33,245.67) -- (565.07,242.63) ;
\draw [color={rgb, 255:red, 80; green, 227; blue, 194 }  ,draw opacity=1 ]   (101.54,260.07) .. controls (89.54,231.07) and (194.33,231.67) .. (322.33,231.67) .. controls (450.33,231.67) and (573.33,215.67) .. (567.33,253.67) ;
\draw    (540.33,259) -- (576.33,249.67) ;
\draw    (524.33,258.67) .. controls (526.33,250.67) and (539.33,242.67) .. (550.33,256.67) ;
\draw  [line width=3] [line join = round][line cap = round] (539.66,259.02) .. controls (539.68,258.02) and (536.07,260.61) .. (537.67,258.31) .. controls (537.93,257.95) and (539.02,257.9) .. (539.01,258.34) ;
\draw    (508.33,83.67) .. controls (510.33,75.67) and (524.33,67.67) .. (533.33,83.67) ;

\draw (526,64) node [anchor=north west][inner sep=0.75pt]   [align=left] {$\displaystyle \pi $};
\draw (170,66) node [anchor=north west][inner sep=0.75pt]   [align=left] {\textcolor[rgb]{0.29,0.56,0.89}{constant area enclosed, length does not increase}};
\draw (101,99) node [anchor=north west][inner sep=0.75pt]   [align=left] {\textcolor[rgb]{0.82,0.01,0.11}{endpoints move in}};
\draw (511,126) node [anchor=north west][inner sep=0.75pt]   [align=left] {tip};
\draw (161,241.33) node [anchor=north west][inner sep=0.75pt]   [align=left] {\textcolor[rgb]{0.29,0.56,0.89}{constant area enclosed, length does not increase}};
\draw (534,305.33) node [anchor=north west][inner sep=0.75pt]   [align=left] {tip};
\draw (513.32,233.82) node [anchor=north west][inner sep=0.75pt]  [font=\normalsize,rotate=-1.77]  {$\sim 179^{\circ }$};
\end{tikzpicture}

   \caption{Endpoint reaching the cone tip (not to scale)}
    \label{fig:enter-label}
\end{figure}
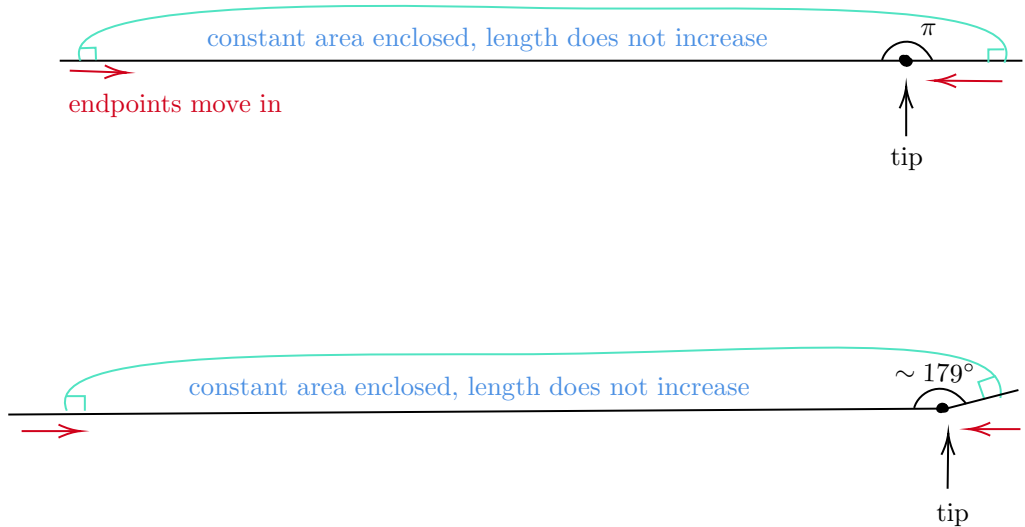

We remark that we do not discount the possibility in this article that an initial curve not satisfying our conditions might be such that under the flow one or both of the endpoints move to the tip of the cone and thus the flow problem ceases to be well-posed.  We suspect such behaviour is possible based on the following example situation.  Consider a `cone' of angle $\pi$, that is, a straight line with an identified `tip'.  Take a very long but rather flat curve with one end quite close to the tip, as indicated in Figure \ref{fig:enter-label}.  The curve diffusion flow with generalised Neumann boundary conditions will keep the area enclosed by the curve and the `cone' constant but decrease the length, so in particular the endpoints will move inward, and the point near the tip could easily reach the tip.  Now, consider a similar situation with cone angle very slightly less than $\pi$.  It seems likely the same behaviour will occur, that is, an end will reach the cone tip, where the boundary condition is not well defined.  

The structure of this article is as follows.  In Section \ref{S:prelim} we set up notation, fully describe our setting and state our main theorems.  In Section \ref{S:flowprops} we provide the evolution equations for various geometric quantities associated with the flow and describe some of the elementary properties of the flow and previous results which are enough to prove the first main theorem.  In Section \ref{S:Kosc} we establish a new key estimate on the oscillation of curvature that then facilitates the proof of our second main theorem by essentially standard arguments which we outline.  We conclude the article with comments concerning the cases of sixth and higher order curve diffusion flows in this setting.

\section{Preliminaries} \label{S:prelim}
\newtheorem{main1}{Theorem}[section]
\newtheorem{main2}[main1]{Theorem}
\newtheorem{PSW2}[main1]{Proposition}
\newtheorem{PSWinfty}[main1]{Proposition}
\newtheorem{interp}[main1]{Proposition}

We consider smooth curves $\alpha$ moving with normal velocity equal to the second arclength derivative of curvature, 
\begin{equation}\label{eq11}
    \frac{\partial}{\partial t} \alpha = k_{ss} \nu \mbox{.}
\end{equation}
Here $\nu$ is the outer unit normal vector field to the curve $\alpha$, and $k=\left<\kappa, \nu  \right>=-\left<  \alpha_{ss}, \nu \right >$  is the scalar curvature of $\alpha$, where $\left<  \cdot, \cdot \right >$ is the usual Euclidean inner product in $\mathbb{R}^2$ and $s$ is the arclength parameter of curve $\alpha$.  The sign in \eqref{eq11} is chosen such that the equation is parabolic in the generalised sense.  Let us further denote by $\tau =\displaystyle{\frac{\alpha_{u} }{\left|\alpha _{u} \right|}}=\alpha _{s}$ the unit tangent vector field along $\alpha$.  Here $u$ is the spatial parameter along the curve, independent of $t$.\\

The boundary conditions for our evolution of planar curves with boundary comprise the classical Neumann condition and a no-curvature-flux condition at each end of the evolving curve.  Fix $0\leq \theta_2 < \theta_1 < \pi$ and let $\gamma_{i} :\left [ 0,\infty  \right )\to \mathbb{R}^{2}~ (i=1,2)$, such that $ \gamma_{i}\left ( \rho  \right )=\rho \left( \cos \theta _{i}, \sin \theta_i \right)$ and denote the image sets by 
$$\bar{\gamma_{i}} :=\left\{\left ( x,y \right )\in \mathbb{R}^{2}: \left ( x,y \right )=\left ( \rho \cos\theta _{i}, \rho\sin \theta _{i}\right ), \rho > 0 \right\}.$$
The set $\bar{\gamma_{1}} \cup \bar{\gamma_{2}}$ forms the boundary for our evolving curve, namely a cone in $\mathbb{R}^2$ with tip at the origin and cone angle $\theta_1 - \theta_2 < \pi$.  

The interior of our initial curve $\alpha_0: \left[ -1, 1\right] \rightarrow \mathbb{R}^2$ is assumed to be contained in the region 
$$\left\{ \left( x, y\right) \in \mathbb{R}^2: \left( x, y\right) = \left( \rho \cos \theta, \rho \sin \theta\right), \rho >0, \theta_2 < \theta < \theta_1 \right\} \mbox{,}$$
 that is, `inside' the cone and such that $\alpha_{0} \left ( -1 \right )\in \bar{\gamma _{1}}$ and $ \alpha_{0} \left ( 1 \right )\in \bar{\gamma _{2}}$.  Moreover, $\alpha_0$ is assumed to meet the boundary at each endpoint perpendicularly (the classical Neumann boundary condition) and be such that the curvature $k_0$ of $\alpha_0$, computed using one-sided derivatives, satisfies
$$\left( k_0\right)_s\left( -1\right) = \left( k_0\right)_s\left( 1 \right) = 0 \mbox{.}$$
This latter condition is the \emph{no curvature flux condition} on the boundary.  Notice that the cone tip itself is not included in the specification of the boundary condition, as the normal to the boundary there is not well defined.  We now let the curve evolve under \eqref{eq11} with \emph{generalised Neumann boundary conditions} as long as this is a well-posed problem.  This means that under the evolution, the classical Neumann boundary condition and no curvature flux condition continue to hold, at least for a short time, by continuity of solutions; the solution will not immediately jump to the cone tip.  It is indeed a key feature of our analysis that we have to avoid the ends of the evolving curve reaching the cone tip, unless they do so in a limiting sense towards the maximal time $T$ of existence of the solution.  Specifically,  $\alpha :\left [ -1,1 \right ]\times \left [ 0,T \right )\to \mathbb{R}^{2}$ is a family of smooth immersed  regular plane curves satisfying \eqref{eq11} with initial curve $\alpha\left ( \cdot ,0 \right )=\alpha _{0}$.  Denoting by $e_-$ and $e_+$ unit vectors perpendicular to $\bar{\gamma_{1}}$ and $\bar{\gamma_{2}}$, the boundary conditions ensure that the ends of the evolving curve continue to meet either side of the cone that hold as long as the solution to \eqref{eq11} exists and
 \begin{equation} \label{E:NBC1}
   \left< \nu\left( -1, t\right), e_- \right> = \left< \nu\left( 1, t\right), e_+ \right> = 0
   \end{equation}
   and
   \begin{equation} \label{E:NBC2}
   k_s\left( \pm 1, t\right) = 0 \mbox{.}
   \end{equation}

\begin{figure}
\centering
\tikzset{every picture/.style={line width=0.75pt}}    
\begin{tikzpicture}[x=0.75pt,y=0.75pt,yscale=-1,xscale=1]

\draw    (54.33,64) -- (200,213) ;

\draw    (334.33,52) -- (200,213) ;

\draw [color={rgb, 255:red, 110; green, 168; blue, 236 }  ,draw opacity=1 ]   (127.17,138.5) .. controls (174.33,108) and (164.33,19) .. (267.17,132.5) ;

\draw [color={rgb, 255:red, 6; green, 82; blue, 173 }  ,draw opacity=1 ][fill={rgb, 255:red, 12; green, 29; blue, 250 }  ,fill opacity=1 ]   (267.17,132.5) -- (294.95,161.55) ;
\draw [shift={(296.33,163)}, rotate = 226.28] [color={rgb, 255:red, 6; green, 82; blue, 173 }  ,draw opacity=1 ][line width=0.75]    (10.93,-3.29) .. controls (6.95,-1.4) and (3.31,-0.3) .. (0,0) .. controls (3.31,0.3) and (6.95,1.4) .. (10.93,3.29)   ;

\draw [color={rgb, 255:red, 7; green, 87; blue, 180 }  ,draw opacity=1 ]   (267.17,132.5) -- (294.07,99.55) ;
\draw [shift={(295.33,98)}, rotate = 129.23] [color={rgb, 255:red, 7; green, 87; blue, 180 }  ,draw opacity=1 ][line width=0.75]    (10.93,-3.29) .. controls (6.95,-1.4) and (3.31,-0.3) .. (0,0) .. controls (3.31,0.3) and (6.95,1.4) .. (10.93,3.29)   ;

\draw    (252.33,133) -- (260.33,140) ;

\draw    (252.33,133) -- (260.33,125) ;

\draw    (136.33,147) -- (144.33,140) ;

\draw    (136.33,131) -- (144.33,140) ;

\draw [color={rgb, 255:red, 74; green, 144; blue, 226 }  ,draw opacity=1 ]   (238.33,214) .. controls (224.16,198.88) and (219.81,200.73) .. (209.24,203.51) ;
\draw [shift={(207.33,204)}, rotate = 345.96] [color={rgb, 255:red, 74; green, 144; blue, 226 }  ,draw opacity=1 ][line width=0.75]    (10.93,-3.29) .. controls (6.95,-1.4) and (3.31,-0.3) .. (0,0) .. controls (3.31,0.3) and (6.95,1.4) .. (10.93,3.29)   ;
 
\draw [color={rgb, 255:red, 74; green, 144; blue, 226 }  ,draw opacity=1 ]   (255.33,213) .. controls (219.63,162.82) and (194.17,177.83) .. (177.15,185.23) ;
\draw [shift={(175.33,186)}, rotate = 337.62] [color={rgb, 255:red, 74; green, 144; blue, 226 }  ,draw opacity=1 ][line width=0.75]    (10.93,-3.29) .. controls (6.95,-1.4) and (3.31,-0.3) .. (0,0) .. controls (3.31,0.3) and (6.95,1.4) .. (10.93,3.29)   ;

\draw (340,30.4) node [anchor=north west][inner sep=0.75pt]    {$\overline{\gamma }_{_{2}}$};

\draw (301,143.4) node [anchor=north west][inner sep=0.75pt]  [color={rgb, 255:red, 35; green, 93; blue, 162 }  ,opacity=1 ]  {$e_{+}$};

\draw (297,96.4) node [anchor=north west][inner sep=0.75pt]  [color={rgb, 255:red, 19; green, 82; blue, 155 }  ,opacity=1 ]  {$\nu $};

\draw (191,48.4) node [anchor=north west][inner sep=0.75pt]  [font=\Large,color={rgb, 255:red, 74; green, 146; blue, 226 }  ,opacity=1 ]  {$\alpha $};

\draw (52,207.4) node [anchor=north west][inner sep=0.75pt]    {$-----------------------$};

\draw (240,173.4) node [anchor=north west][inner sep=0.75pt]    {$\theta _{1}{}$};

\draw (229.07,192.33) node [anchor=north west][inner sep=0.75pt]  [rotate=-0.45]  {$\theta _{2}$};

\draw (42,35.4) node [anchor=north west][inner sep=0.75pt]   {$\overline{\gamma }_{_{1}}$};
\end{tikzpicture}
\caption{The set-up}
\end{figure}
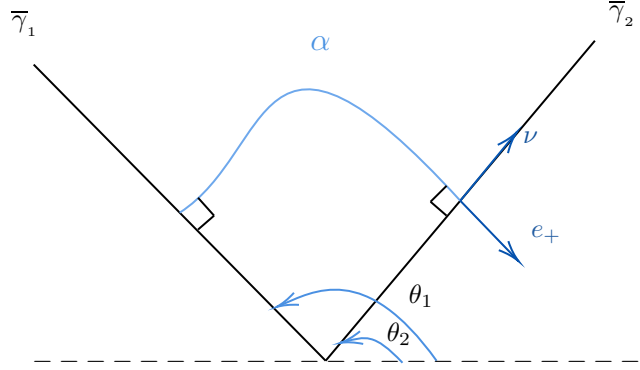

 It is convenient in the analysis to consider the evolving curves parametrised by their arc length $s$, which itself is a function of $t$.  The arclength parameter $s$ is related to parameter $u\in\left[ -1, 1 \right]$ by
\begin{equation} \label{E:s}
  s\left( u\right) = \int_{-1}^{u} \left| \alpha_u\left( \tilde u\right) \right| d\tilde u \mbox{.}
  \end{equation}
 This means that we measure the arc length from the point $\alpha\left( -1\right)$ on curve $\alpha$ so, in particular, $s\left( -1\right) = 0$.  We further have, from \eqref{E:s}, that $ds=\left| \alpha _{u}\right|du$.  By slight abuse of notation, we will denote quantities associated with the evolving curve using the same symbols, be they functions of $u$ and $t$, or $s$ and $t$, as is customary, and should not lead to any confusion.  When we evaluate at endpoints we will denote the spatial argument using $\pm 1$ but spatial derivatives will be respect to $s$ and interpreted in the appropriate one-sided sense.  We denote by $k_{s^{n}}$ the $n$-th iterated derivative of $k$ with respect to arclength and write $k_{s^{n}}^{2}$ for $\left ( k_{s^{n}} \right )^{2}$. 
 
 Throughout the article, we use $L$ to denote the length of $\alpha\left( \cdot, t\right)$ and write
\begin{equation}\label{eql}
    L\left( t\right) = L[\alpha\left( \cdot, t\right) ]=\int\, ds.
\end{equation}
All integrals will be over the curve $\alpha$ unless otherwise indicated.  The region bounded by the curve and the cone is an area bounded by a piecewise smooth closed curve.  This area may be computed using Green's Theorem (which can be considered a special case of both Stokes' Theorem and the Divergence Theorem).  On the rays of the cone, the position vector is perpendicular to the normal vector, so the result is that only the line integral along $\alpha$ is nonzero and we have
\begin{equation} \label{E:area}
  A\left( t\right) = A\left[\alpha\left( \cdot, t\right) \right] = +\frac{1}{2} \int \left< \alpha, \nu \right> ds \mbox{,}
\end{equation}
where the sign is $+$ due to our choice of outer unit normal.\\

\noindent \emph{Remark:} As an alternative to the use of Green's Theorem above, if the curve $\alpha$ is star-shaped, then \eqref{E:area} results from the use of polar coordinates.  The use of Green's Theorem is preferred as we do not wish to restrict $\alpha$ to being star-shaped.\\

The average curvature of $\alpha\left( \cdot, t\right)$  is defined as\begin{equation}\label{eqkk}
 \bar{k}[\alpha]=\displaystyle{\frac{1}{L}}\int k\,ds.   
\end{equation}

Additionally, the oscillation of curvature is given by
\begin{equation}\label{eqk}
K_{\mbox{osc}}[\alpha ]=L\int \left ( k-\bar{k} \right )^{2}ds 
\end{equation}
and the rotation number of $\alpha\left( \cdot, t\right)$ is
\begin{equation} \label{E:omega}
  \omega \left [ \alpha  \right ] :=\frac{1}{2\pi }\int k\,ds.
  \end{equation}

Equation \eqref{E:omega} may be thought of as the general definition of rotation number for any curve.  As the rotation number gives the net turning of the tangent vector to $\alpha$ we observe that in our particular setting 
$$\omega = \frac{\theta_1 - \theta_2}{2\pi} \mbox{.}$$
In view of our constraints on $\theta_1$ and $\theta_2$, we have $\omega \in \left( 0, \frac{1}{2} \right)$.\\

 With the above setting described we can now state our main results.
 
 \begin{main1} \label{T:main1}
 Let $\alpha_0: \left[ -1, 1\right]\to\mathbb{R}^{2}$ be a smooth curve of length $L_0>0$, contained within the cone as described above, with rotation number $\omega \in \left( 0, \frac{1}{2} \right)$ and compatible with the generalised Neumann boundary conditions \eqref{E:NBC1} and \eqref{E:NBC2}.  There is a unique solution $\alpha\left( \cdot, t\right)$ to the curve diffusion flow, \eqref{eq11} with generalised Neumann boundary conditions \eqref{E:NBC1} and \eqref{E:NBC2} and initial condition $\alpha\left( \cdot, 0 \right) = \alpha_0$ on a maximal time interval $\left[ 0, T\right)$, $T\leq \infty$.  
 
 If $T< \infty$ then at least one of the following is true:
 \begin{itemize}
   \item at least one end of the evolving curve $\alpha\left( \cdot, t\right)$ has reached the cone tip;
   \item the curvature of the curve satisfies
$$\int k^2 ds \geq c\left( T-t\right)^{-\frac{1}{4}}$$
 for some constant $c$.\\
 \end{itemize}
 \end{main1}
 
 \noindent \emph{Remark:} The possibility of $T=\infty$ is addressed in Theorem \ref{T:main2}, however, it is also clear that any circular arc is a stationary solution to \eqref{eq11} satisfying the generalised Neumann boundary conditions.  Moreover a standard stability argument shows that any suitable initial curve that is also $C^{4, \alpha}$-close to a circular arc and suitably far from the tip of the cone will provide a solution that exists for all time and converges exponentially to a circular arc in the $C^\infty$-topology.  The next result weakens the closeness requirement to a smallness condition on the oscillation of curvature of the initial curve.\\
 
 \begin{main2} \label{T:main2}
 Let $\alpha_0: \left[ -1, 1\right]\to \mathbb{R}^{2}$ be as in the statement of Theorem \ref{T:main1} but also satisfying
\begin{equation} \label{E:smallness}
  K_{\mbox{osc}} \leq  \left ( \frac{\pi }{12} \right )^{2}\left[ -24\,\omega +\sqrt{\left(24\, \omega\right) ^{2}+\frac{48}{\pi }\left( 1- 4\, \omega ^{2} \right) }\right ]^{2} 
\end{equation}
with ends sufficiently far from the tip of the cone.
Then, \eqref{eq11} with the generalised Neumann boundary conditions \eqref{E:NBC1} and \eqref{E:NBC2} and initial condition $\alpha\left( \cdot, 0 \right) = \alpha_0$ has a unique solution $\alpha\left( \cdot, t\right)$ for all time.  Condition \eqref{E:smallness} is preserved under the flow, $\alpha\left( \cdot, t\right)$ is smooth for all time and converges exponentially in the $C^\infty$-topology as $t\rightarrow \infty$ to a circular arc that, together with the cone, bounds the same area as that bounded by the initial curve and the cone.
\end{main2}

\noindent \emph{Remark:} That the solution converges in the $C^\infty$-topology to the circular arc denoted $\alpha_\infty: \left[ -1, 1\right] \rightarrow \mathbb{R}^2$ means that
$$\lim_{t\rightarrow \infty} \left| \alpha\left( u, t\right) - \alpha_\infty\left( u \right) \right| = 0$$
uniformly in $u$ and for all $\ell \in \mathbb{N},$
$$\lim_{t\rightarrow \infty} \left| \alpha_{u^\ell}\left( u, t\right) \right| = 0$$
uniformly in $u$.  That the convergence is exponential means there exist absolute constants $C_\ell, \delta_\ell >0$, $\ell\in \mathbb{N} \cup \left\{ 0 \right\}$ such that
$$\left| \alpha\left( u, t\right) - \alpha_\infty\left( u \right) \right| \leq C_0 \, e^{-\delta_0, t}$$
and
$$\left| \alpha_{u^\ell}\left( u, t\right) \right| \leq C_\ell \, e^{-\delta_\ell\, t} \mbox{.}$$

We note that typically exponential convergence of higher derivatives is obtained via interpolation and the constants $\delta_\ell$ are not all the same.  To obtain exponential convergence with the same constants one may use a linearisation argument.\\


To prove our main theorems we will require some fundamental tools that are all well known.  We will state them here with relevant references for the convenience of the reader.  First are the standard Poincaré-Sobolev-Wirtinger [PSW] inequalities for  curves with boundary. 

\begin{PSW2} \label{psw}
Let $L>0$. Let $g:[0,L]\rightarrow \mathbb{R}$ be an absolutely continuous function with $\int_{0}^{L}g(x) \,dx=0.$ Then 
\begin{equation}\label{eqpsw}
    \int_{0}^{L}g^{2}(x)\, dx\leq \frac{L^{2}}{\pi ^{2}}\int_{0}^{L}\left( g_{x} \right)^{2}dx.
\end{equation}
Similarly, if $g:[0,L]\rightarrow \mathbb{R}$ is absolutely continuous and $g(0)=g(L)=0$ then 
\begin{equation*}
    \int_{0}^{L}g^{2}(x) \,dx\leq \frac{L^{2}}{\pi ^{2}}\int_{0}^{L}\left( g_{x} \right)^{2}(x)dx.
\end{equation*}
\end{PSW2}

\begin{PSWinfty} \label{eq266}
If $g:[0,L]\rightarrow \mathbb{R}$ is absolutely continuous and $g(0)=g(L)=0$ then 
\begin{equation}\label{pp}
\left \| g \right \|_{\infty }^{2}\leq \frac{L}{\pi }\int_{0}^{L}\left( g_{x} \right)^2dx.
\end{equation}
Similarly, if  $g:[0,L]\rightarrow \mathbb{R}$ is absolutely continuous and $\int_{0}^{L}g dx=0$ then
\begin{equation}\label{eqppsw}
\left \| g \right \|_{\infty }^{2}\leq \frac{2L}{\pi }\int_{0}^{L}\left( g_{x} \right)^2 dx.
\end{equation}
\end{PSWinfty}
For the next interpolation inequality we require, we first set up
 some notation.  For $i_1, \ldots, i_n  \in \mathbb{N}\cup \left\{ 0\right\}$ such that $i_1 + \ldots+ i_n=m$ , denote 
 by $P_n^m\left( k\right)$ any linear combination of terms of type $k_{s^{i_1}} k_{s^{i_2}} \cdots k_{s^{i_n}}$.  Here we are making the identification $k_{s^0}= k$.
 
  It is convenient to use the following scale-invariant norms: we define
 $$\left\| k \right\|_{\ell, p} := \sum_{i=0}^\ell \left\| k_{s^i} \right\|_p$$
 where
 $$\left\| k_{s^i} \right\|_p := L^{i+1-\frac{1}{p}} \left( \int \left| k_{s^i} \right|^p ds\right)^{\frac{1}{p}} \mbox{.}$$
 
 The following interpolation inequality for curves with boundary appears as
 \cite{DP14}*{Lemma 4.3}.

\begin{interp}\label{p}
  Let $\alpha:\left [ -1, 1 \right ]\to \mathbb{R}^{2}$  be a smooth curve with boundary. Then for any term $P_{n}^{m}\left ( k \right )$ with $n\geq 2$  that contains derivatives of $k$ of order at most $\ell-1$,
  $$\int \left|P_{n}^{m}\left ( k \right ) \right| ds\leq c\, L^{\ell-m-n}\left\| k\right\|_{0,2}^{n-p}\left\| k\right\|_{\ell,2}^{p},$$
  where $p=\displaystyle{\frac{1}{\ell}\left ( m+\frac{1}{2} n-1\right )}$ and $c=\left ( \ell,m,n \right ) $. Furthermore, if $m+\frac{n}{2}< 2\ell+1$ then $p<2$ and for any $\varepsilon>0,$ there is another constant $c=c\left( n, \ell, m \right)$ such that
  $$\int\left|P_{n}^{m}\left ( k \right ) \right| ds\leq \varepsilon \int  k_{s^{\ell}}^{2}ds+c\,\varepsilon^{\frac{-p}{2-p}}\left( \int k^{2}ds  \right)^{\frac{n-p}{2-p}}+c \left(\int k^{2}ds  \right)^{m+n-1}.$$
\end{interp}

  \section{Preliminary properties of the flow} \label{S:flowprops}
  \newtheorem{evlneqns1}{Lemma}[section]
  \newtheorem{BCs}[evlneqns1]{Lemma}
  \newtheorem{evlneqns2}[evlneqns1]{Corollary}
  \newtheorem{rot}[evlneqns1]{Corollary}
  \newtheorem{blowup}[evlneqns1]{Theorem}
  \newtheorem{Area}[evlneqns1]{Corollary}
  \newtheorem{Lbound}[evlneqns1]{Corollary}
  
  We begin this section with evolution equations for geometric quantities needed in our analysis.  These are derived in a straightforward way, as in \cites{W13, WW} for example.

\begin{evlneqns1} \label{T:evlneqns1}
Under the flow \eqref{eq11} we have the following evolution equations for various pointwise geometric quantities.
\begin{itemize}
  \item[\textnormal{(i)}] $\frac{\partial}{\partial t} ds = k\, k_{ss} ds$;
  \item[\textnormal{(ii)}] $\frac{\partial}{\partial t} k_{s^\ell} = -k_{s^{\ell+4}} - \sum_{j=0}^{\ell} \left( k\, k_{s^{\ell-j}} k_{ss} \right)$ for each $\ell = 0, 1, 2, \ldots$.
  
  \end{itemize}
  \end{evlneqns1}
  
  Exactly as in \cite{WW} we have that all odd curvature derivatives are equal to zero on the boundary, as long as the solution to \eqref{eq11} exists.  This is critical for future calculations where we `integrate by parts', that is, apply the Divergence Theorem to various functions involving curvature and its derivatives on our curve with boundary.
  
  \begin{BCs} \label{T:BCs}
  Under the flow \eqref{eq11} with generalised Neumann boundary conditions as stated, we have on the boundary for each $\ell \in \mathbb{N}$ and all $t\in \left[ 0, T\right)$
  $$k_{s^{2\ell-1}}\left( \pm 1, t\right)=0 \mbox{.}$$
  \end{BCs}
  
  Next we state the evolution equations for some integral geometric quantities.  Again these are straightforward to derive from the definitions and Lemma \ref{T:evlneqns1}.  
  
  \begin{evlneqns2} \label{T:evlneqns2}
  Under the flow \eqref{eq11} we have the following evolution equations for various integral geometric quantities.
\begin{itemize}
  \item[\textnormal{(i)}] $\frac{\mathrm{d} }{\mathrm{d} t} L = - \int  k_s^2 ds$;
  \item[\textnormal{(ii)}] $\frac{\mathrm{d} }{\mathrm{d} t}A = 0$;
  \item[\textnormal{(iii)}] $\frac{\mathrm{d} }{\mathrm{d} t} \overline{k} = \frac{2\, \pi\, \omega}{L} \int k_s^2 ds$;
  \item[\textnormal{(iv)}] $\frac{\mathrm{d} }{\mathrm{d} t}\int k^2 ds = - \int k_{ss}^2 ds + 3 \int k^2 k_s^2 ds$;
  \item[\textnormal{(v)}] $\frac{\mathrm{d} }{\mathrm{d} t} K_{\mbox{osc}} = - K_{\mbox{osc}} \frac{\left\| k_s \right\|_2^2}{L} - 2\, L\, \left\| k_{ss} \right\|_2^2 + 3\, L\, \int\left( k - \overline{k}\right)^2 k_s^2 ds+ 6\, \overline{k} \, L\, \int \left( k - \overline{k} \right) k_s^2 ds$
  
  \hspace*{\fill}$ + 2\, \overline{k}^2 L \, \left\| k_s \right\|_2^2$.
  \end{itemize}
  \end{evlneqns2}
  

  
  We complete this section with some fundamental properties of the curve diffusion flow in our setting which are consistent with other settings.  For these two results we do not require a smallness condition on the initial oscillation of curvature.  The first result is that the rotation number $\omega$ is constant under the flow.  
  
  \begin{rot} \label{T:omega}
 Under the flow \eqref{eq11} with generalised Neumann boundary conditions \eqref{E:NBC1} and \eqref{E:NBC2} we have
 \begin{equation}
     \frac{\mathrm{d} }{\mathrm{d} t}\int _{\alpha }k\, ds=0.
 \end{equation}
\end{rot}

\noindent \textbf{Proof:} We use Lemma \ref{T:evlneqns1} and integration by parts, noting Lemma \ref{T:BCs}, to calculate
 \begin{multline*}
\frac{\mathrm{d} }{\mathrm{d} t}\int _{\alpha }k\, ds =\int _{\alpha }\frac{\partial }{\partial t}k\, ds+\int _{\alpha }k\frac{\partial }{\partial t}ds 
 =\int _{\alpha }\left ( -F_{ss}-k^{2}F \right )ds+\int _{\alpha }k^{2}Fds\\
 = -\int _{\alpha }F_{ss}ds= -F_{s}\bigg|_{\left\{ 0,L\right\}}= 0.
\end{multline*}
In view of the definition \eqref{E:omega} the proof is complete.\hspace*{\fill}$\Box$\\

\noindent \emph{Remark:} The result of Corollary \ref{T:omega} is of course expected in view of the earlier discussion of the rotation number.  However, we have included it to highlight two points.  The first is the role of the boundary conditions in the above calculation.  The second is to confirm that the rotation number will remain constant while the smooth solution to the flow exists.\\ 

In the case of closed curves evolving by the curve diffusion flow, the area bounded by the evolving curve is constant \cite{W13}[Lemma 3.1].  A similar result holds here: the area bounded by the curve and the two ray segments emanating from the cone tip is constant. 

\begin{Area} \label{T:area}
Under the flow \eqref{eq11} with generalised Neumann boundary conditions \eqref{E:NBC1} and \eqref{E:NBC2}, the area enclosed by the curve and the cone is constant.
\end{Area}

\noindent \textbf{Proof:} This follows by the same argument as in the proof of \cite{W13}[Lemma 3.1] using Lemma \ref{T:evlneqns2} (ii).  The only difference in the calculation is the appearance of boundary terms in integrating by parts, but these are equal to zero in view of Lemma \ref{T:BCs}.\hspace*{\fill}$\Box$\\

\begin{Lbound} \label{T:Lbound}
Under the flow \eqref{eq11} with generalised Neumann boundary conditions \eqref{E:NBC1} and \eqref{E:NBC2}, the length of the evolving curve is uniformly bounded above and below by positive constants.
\end{Lbound}

\noindent \textbf{Proof:} As \eqref{eq11} is the $H^{-1}$-gradient flow for the length functional, we have directly from Lemma \ref{T:evlneqns2}, (i) that
$$L\left( t\right) \leq L_0 \mbox{.}$$
A positive lower bound $\underline{L}$ on $L\left( t\right)$ follows by contradiction: if $L\left( t\right)$ is too small, then the region bounded by $\alpha\left( \cdot, t\right)$ and the two rays of the cone could not enclose area $A_0$, which is constant under the flow.  Hence there exists $\underline{L}>0$ such that 
$$L\left( t\right) \geq \underline{L} > 0$$
holds on $\left[ 0, T\right)$. \hspace*{\fill}$\Box$\\

\noindent \textbf{Proof of Theorem \ref{T:main1}:} 
The case of at least one end of the evolving curve approaching the tip is clear, as the problem at such a time ceases to be well-posed.  Suppose $T$ is finite but both ends of the evolving curve remain away from the tip of the cone.  Applying Proposition \ref{p} to the second term on the right hand side of Corollary \ref{T:evlneqns2}, (iv) we obtain a blow-up criterion for the flow, similarly as in the case of closed curves evolving by curve diffusion \cite{DKS02}[Theorem 3.1] and as stated above.\hspace*{\fill}$\Box$\\


\section{Controlling the oscillation of curvature} \label{S:Kosc}
\newtheorem{Koscmonotone}{Proposition}[section]
\newtheorem{conepoint}[Koscmonotone]{Corollary}
\newtheorem{Lbound2}[Koscmonotone]{Corollary}
\newtheorem{L2bounds}[Koscmonotone]{Corollary}
\newtheorem{expdecay}[Koscmonotone]{Proposition}

We begin by showing that if $\alpha_0$ has sufficiently small oscillation of curvature, then this continues to hold under \eqref{eq11}.  The calculation is similar to that of \cite{WW} however here we do not have $\omega=0$ and hence neither do we have $\overline{k}=0$.  Thus we have to deal with some additional terms in the evolution equation for $K_{\mbox{osc}}$.

\begin{Koscmonotone} \label{T:Koscmonotone}
 Let  $\alpha :\left [ -1,1 \right ]\times \left [ 0,T \right )\to \mathbb{R}^{2}$ be a solution to \eqref{eq11} with $\omega < \frac{1}{2}$, generalised Neumann boundary conditions as described earlier and initial curve $\alpha_0$ satisfying \eqref{E:smallness}.  Then \eqref{E:smallness} remains true for $t\in \left[ 0, T\right)$.
 \end{Koscmonotone}
   
\noindent \textbf{Proof:}  We estimate the terms on the right hand side of Corollary \ref{T:evlneqns2} (v).  Using Proposition \ref{eq266}  and the H\"{o}lder inequality we obtain
$$3\,L\int _{\alpha }\left ( k-\bar{k} \right )^{2}k_{s}^{2}ds\leq 3\, L \left\|k_{s} \right\|_{\infty }^{2}\int _{\alpha }\left ( k-\bar{k} \right )^{2}ds\leq \frac{6\,L}{\pi }K_{\mbox{osc}}\left\|k_{ss} \right\|_{2}^{2},$$
 $$6\,\bar{k}\,L\int _{\alpha }\left ( k-\bar{k} \right )k_{s}^{2}ds\leq 12\,\omega \,\pi \sqrt{K_{\mbox{osc}}}\left\|k_{s} \right\|_{\infty }^{2}\leq 24\, \omega\, L \sqrt{K_{\mbox{osc}}} \left\|k_{ss} \right\|_{2 }^{2},$$
 and
$$ 2\, L\, \bar{k}^{2}\left\|k_{s} \right\|_{2}^{2}=2\, L\, \left ( \frac{2\, \omega \pi }{L} \right )^{2}\left\|k_{s} \right\|_{2}^{2} \leq 8\, \omega ^{2}L\left\|k_{ss} \right\|_{2}^{2}.$$
Substituting these into Corollary \ref{T:evlneqns2} (iv) we obtain
  \begin{equation} \label{E:Kosc}
    \frac{\mathrm{d} }{\mathrm{d} t}K_{\mbox{osc}}+K_{\mbox{osc}}\frac{\left\|k_{s} \right\|_{2}^{2}}{L}+\left [2-\frac{6}{\pi }K_{\mbox{osc}}-24\, \omega \sqrt{K_{\mbox{osc}}}-8\,\omega ^{2}\right ]L\left\|k_{ss} \right\|_{2}^{2}
 \leq 0.
 \end{equation}

Given $\omega < \frac{1}{2}$, the result follows for $K_{\mbox{osc}}$ satisfyinging \eqref{E:smallness}.\hspace*{\fill}$\Box$\\

Proposition \ref{T:Koscmonotone} ensures that unless one or both end reach the tip of the cone, the solution to \eqref{eq11} with generalised Neumann boundary conditions \eqref{E:NBC1} and \eqref{E:NBC2} exists for all time.

\begin{conepoint} \label{T:conepoint}
Suppose $\alpha :\left [ -1,1 \right ]\times \left [ 0,T \right )\to \mathbb{R}^{2}$ is a solution to \eqref{eq11} with $\omega < \frac{1}{2}$, generalised Neumann boundary conditions \eqref{E:NBC1} and \eqref{E:NBC2} and initial curve $\alpha_0$ satisfying \eqref{E:smallness}.  Then unless an end of the curve reaches the cone tip, $T=\infty$.
\end{conepoint}

\noindent \textbf{Proof:} Suppose for the sake of contradiction that $T<\infty$ and the ends of the curve remain away from the cone tip.  Theorem \ref{T:main1} gives that $\int k^2 ds$ becomes unbounded.  However from Proposition \ref{T:Koscmonotone} we also have
$$ L \int \left( k - \overline{k}\right)^2 ds = L \int k^2 ds - 4\pi^2 \omega^2 = K_{\mbox{osc}} \left( t\right) \leq K_{\mbox{osc}}\left( 0\right)$$
and therefore
\begin{equation} \label{E:k2bound}
  \int k^2 ds \leq \frac{K_{\mbox{osc}}\left( 0\right)+4\pi^2 \omega^2}{L}  \mbox{.}
\end{equation}
In view of this inequality, the only way $\int k^2 ds$ can become unbounded as $t\rightarrow T$ is if $\lim_{t\rightarrow T} L\left( t\right) = 0$.  This, however, is impossible, by the uniform positive lower length bound of Corollary \ref{T:Lbound}. \hspace*{\fill}$\Box$\\

%

In view of Corollary \ref{T:Lbound}, \eqref{E:k2bound} shows that $\int k^2 ds \leq C_0^2$ for all $t$.  Using Proposition \ref{p} similarly as in  \cite{DP14} implies all curvature derivatives are bounded in $L^2$.

\begin{L2bounds} \label{T:L2bounds}
Suppose $\alpha :\left [ -1,1 \right ]\times \left [ 0,\infty \right)\to \mathbb{R}^{2}$ is a solution to \eqref{eq11} with $\omega < \frac{1}{2}$, generalised Neumann boundary conditions as described earlier and initial curve $\alpha_0$ satisfying \eqref{E:smallness}.  For each $m\in \mathbb{N}\cup \left\{ 0\right\}$ there exists a constant $C_m^2 >0$ such that
$$\int k_{s^m}^2 ds \leq C_m^2$$
for all $t\in\left[ 0, \infty \right)$.
\end{L2bounds}

Since $\omega < \frac{1}{2}$, a refinement in the proof of Proposition \ref{T:Koscmonotone} yields a crucial exponential decaying quantity under the flow, namely $\int \left( k - \overline{k}\right)^2 ds$.

\begin{expdecay} \label{T:expdecay}
Suppose $\alpha :\left [ -1,1 \right ]\times \left [ 0,\infty \right)\to \mathbb{R}^{2}$ is a solution to \eqref{eq11} with $\omega < \frac{1}{2}$, generalised Neumann boundary conditions as described earlier and initial curve $\alpha_0$ satisfying \eqref{E:smallness}.  Then
$$\int\left( k-\overline{k}\right)^2 ds =\frac{K_{\mbox{osc}}\left( t\right)}{L\left( t\right)} \leq \frac{K_{\mbox{osc}}\left( 0\right)}{L_0} e^{- \frac{\delta \pi^4}{L_0^4} t} \mbox{.}$$
\end{expdecay}

\noindent \textbf{Proof:} Under the stated conditions we know that $K_{\mbox{osc}}$ does not increase under the flow so
$$2 - \frac{6}{\pi} K_{\mbox{osc}} - 24\omega \sqrt{K_{\mbox{osc}}} \geq \delta >0$$
continues to hold and from \eqref{E:Kosc} we obtain
 \begin{equation*} 
    \frac{\mathrm{d} }{\mathrm{d} t}K_{\mbox{osc}}+K_{\mbox{osc}}\frac{\left\|k_{s} \right\|_{2}^{2}}{L} \leq - \delta L\left\|k_{ss} \right\|_{2}^{2} \mbox{.}
 \end{equation*}
Hence, using Proposition \ref{psw} twice and Corollary \ref{T:Lbound},
$$\frac{d}{dt} K_{\mbox{osc}} + K_{\mbox{osc}} \frac{\left\| k_s \right\|_2^2}{L} \leq - \delta \, L \frac{\pi^4}{L^5} K_{\mbox{osc}} \leq - \delta \frac{\pi^4}{L_0^4} K_{\mbox{osc}}  \mbox{.}$$
Hence, in view of Corollary \ref{T:evlneqns2}, (i),
$$\frac{d}{dt} \ln \frac{K_{\mbox{osc}}}{L} \leq - \delta \, \frac{\pi^4}{L_0^4}$$
from which the result follows.\hspace*{\fill}$\Box$\\

Proposition \ref{T:expdecay} allows us to establish an improved upper bound on $L$ under the flow.  Although we do not actually need this here, we include it as we feel this kind of estimate could be useful in other settings.

\begin{Lbound2} \label{T:Lbound2}
Suppose $\alpha :\left [ -1,1 \right ]\times \left [ 0,T \right )\to \mathbb{R}^{2}$ is a solution to \eqref{eq11} with $\omega < \frac{1}{2}$, generalised Neumann boundary conditions \eqref{E:NBC1} and \eqref{E:NBC2} and initial curve $\alpha_0$ satisfying \eqref{E:smallness}.  Then
$$L^3\left( t\right) \leq L_0^3 \left[ 1 + \frac{3\, K_{\mbox{osc}}\left( 0\right)}{\delta \pi^2} \left( e^{-\frac{\delta \, \pi^4}{L_0^4} t} - 1\right) \right] \mbox{.}$$
\end{Lbound2}

\noindent \textbf{Proof:} Using Corollary \ref{T:evlneqns2}, (i) and Proposition \ref{psw} we estimate
$$ \frac{\mathrm{d} }{\mathrm{d} t} L = - \int k_s^2 ds \leq -\frac{\pi^2}{L^2} \int \left( k - \overline{k}\right)^2 ds \leq -\frac{\pi^2}{L^2} \frac{K_{\mbox{osc}}\left( 0\right)}{L_0} e^{-\frac{\delta \, \pi^4}{L_0^4} t} \mbox{,}$$
where for the last step we used Proposition \ref{T:Koscmonotone}.  Hence
$$ \frac{\mathrm{d} }{\mathrm{d} t} L^3 \leq - 3\, \pi^2 \frac{K_{\mbox{osc}}\left( 0\right)}{L_0} e^{-\frac{\delta \, \pi^4}{L_0^4} t} \mbox{,}$$
from which the result follows. \hspace*{\fill}$\Box$\\

\noindent \textbf{Completion of the proof of Theorem \ref{T:main2}:} Corollary \ref{T:conepoint} gives that $T=\infty$ unless an end of the curve has reached the cone tip.  Assuming then that neither end of the evolving curve reaches the cone tip, we have $T=\infty$.  Using then Corollary \ref{T:L2bounds} and Proposition \ref{T:expdecay} we obtain by a standard induction argument using integration by parts, Lemma \ref{T:BCs} and the H\"{o}lder inequality that all curvature derivatives decay in $L^2$ and in $L^\infty$ to zero.  The first step of this argument is
$$\int k_s^2 ds = - \int \left( k - \overline{k} \right) k_{ss} ds \leq \left[ \int \left( k - \overline{k} \right)^2 ds \right]^{\frac{1}{2}} \left( \int k_{ss}^2 ds \right)^{\frac{1}{2}} \leq C_2 \left( \frac{K_{\mbox{osc}}\left( 0\right)}{L_0} \right)^{\frac{1}{2}} e^{- \frac{\delta \pi^4}{2\, L_0^4} t } \mbox{.}$$
Via Proposition \ref{eq266} inequality one obtains that $k \rightarrow \overline{k}$ pointwise and all derivatives of $k$ approach zero pointwise exponentially.  

For the sake of completeness, we provide the inductive step.  Suppose then $\int k_{s^m}^2 ds \leq \tilde C_m e^{-\delta_m t}$ for some $\tilde C_m, \delta_m >0$.  Integrating by parts and using the boundary condition Lemma \ref{T:BCs}, we estimate using the H\"{o}lder inequality
$$\int k_{s^{m+1}}^2 ds = - \int k_{s^{m+2}} k_{s^m} ds \leq \left( \int k_{s^{m+2}} ds\right)^{\frac{1}{2}} \left( \int k_{s^{m}} ds\right)^{\frac{1}{2}} \mbox{.}$$
Using now the inductive assumption and Corollary 4.3 yields
$$k_{s^{m+1}}^2 ds \leq C_{m+2} \tilde C_m^{\frac{1}{2}} e^{-\frac{\delta_m}{2} t} \mbox{.}$$
It follows that each order of curvature derivative converges to zero exponentially in $L^2$.  The corresponding pointwise convergence follows as usual using Lemma \ref{eq266} and the length bound.  To obtain exponential decay with the same constant one can consider the linearisation about the limiting circular arc.

To complete the proof that the evolving curve $\alpha\left( \cdot, t\right)$ converges exponentially to the unique circular arc $\alpha_\infty\left( \cdot\right)$ with curvature $\overline{k}$, in terms of convergence of the derivatives of $\alpha$, one can follow the argument in \cite{WW}, for example.  The limiting arc $\alpha_\infty$, together with the boundary rays of the cone, encloses the same area as the initial curve and the cone since from Corollary \ref{T:evlneqns2}, (ii), the enclosed area is constant under the flow.  

The first step is to control the length of the tangent vector in the $u$ parametrisation.  In view of \eqref{eq11}, we have
$$\frac{\partial}{\partial t} \left| \alpha_u \right| = k\, k_{ss} \left| \alpha_u \right| \mbox{.}$$
It follows, in view of the curvature bounds, exponential decay and Lemma \ref{psw}, there are constants $C, \delta >0$ such that for each fixed $u\in \left[ -1, 1\right]$ (with, as usual, one-sided derivatives at the endpoints),
$$-C\, e^{-\delta \, t} \leq \frac{\partial}{\partial t} \ln \left| \alpha_u\left( u, t\right) \right| \leq C\, e^{-\delta\, t} \mbox{.}$$
Integrating the above and estimating, there is a constant $c_0>0$ such that
$$\frac{1}{c_0} \left| \alpha_u\left( u, 0\right) \right| \leq \left| \alpha_u\left( u, t\right) \right| \leq c_0 \left| \alpha_u\left( u, 0\right) \right| \mbox{.}$$
This uniform control can then be used to convert control on arclength derivatives to control on parameter derivatives, again appealing to an ODE argument for any fixed $u$.  Specifically, it can be shown by induction, exactly as in \cite{WW}, that for each $\ell \in \mathbb{N}$ and every $\left( u, t\right)$,
$$\left( \left| \alpha_u\right|^{-1} \right)_{u^\ell} \left( u, t\right) \leq c_\ell \mbox{ and } \left| \left( k\, k_{ss} \right)_{u^\ell} \right| \left( u, t\right) \leq c_\ell e^{-\delta\, t} \mbox{.}$$
Exponential decay of $\left\| \alpha_{u^\ell} \right\|_{L^\infty\left( -1, 1\right)}$ for each $\ell \in \mathbb{N} \backslash \left\{ 1\right\}$ now follows using the above and converting $u$ derivatives of $\alpha$ into $s$ derivatives of $k\,\nu$ using the Frenet equations and using the exponential convergence of the curvature derivatives and derivatives of the speed as well.

Finally we note that the distance travelled by any point on the initial curve $\alpha_0$ is finite, as can be shown by integrating the normal speed (as in \cite{MWW22} and in \cite{WW}, for examples).  Specifically we estimate from \eqref{eq11}
$$\frac{\partial}{\partial t} \alpha \leq \left| k_{ss} \nu \right| \leq c\, e^{-\delta_0\, t}$$
and so by integration for any particular $u$,
$$\left| \alpha\left( u, t\right) \right| \leq \left| \alpha_0\left( u \right) \right| + \frac{c}{\delta_0} \mbox{.}$$
This shows that $\left| \alpha\left( u, t\right)\right|$ remains bounded and, moreover, we have a bound on how far any point on the initial curve can move under the flow.  Thus we can say that if $\alpha_0$ is initially far enough from the tip of the cone, neither end will reach the tip under the flow, so we can be sure $T=\infty$.\hspace*{\fill}$\Box$\\

\noindent \emph{Remark:}  One can consider polyharmonic curve diffusion flows in the cone for initial curves of rotation number $\omega < \frac{1}{2}$ as an extension of the situation here.  In these case
$$\frac{\partial}{\partial t} \alpha = \left( -1\right)^{m+1} k_{s^{2m}} \nu \mbox{,}$$
which is a $\left( 2m+2\right)$-th order parabolic evolution equation.  Such an equation was considered in the case of closed curves in \cite{PW16} and for curves with boundary on parallel lines in \cite{MWW22}.  Here we can adopt the same (further) generalised Neumann boundary conditions as in \cite{MWW22}; we specify all odd derivatives of curvature up to order $2m+1$ are equal to zero on the boundary.  Lemma \ref{T:BCs} then again holds for all higher odd derivatives on the boundary as long as the solution exists.  Analogues of Theorem \ref{T:main1} and Theorem \ref{T:main2} hold.  In the case $m=1$ it turns out that exactly the same smallness condition \eqref{E:smallness} is preserved under the flow.  The one additional term in the evolution equation for $K_{\mbox{osc}}$ in this case has a good sign, so causes no harm.  In the cases of higher $m$ there are more terms in the evolution equation for $K_{\mbox{osc}}$, however, a result of the form \cite{PW16}[Lemma 7] holds and so for $\omega < \frac{1}{2}$ a similar result again holds at least if $K_{\mbox{osc}}\left( 0\right)$ is initially sufficiently small.

\end{document}